\input amstex
\documentstyle{amsppt}
\magnification=\magstep0
\define\cc{\Bbb C}
 \define\z{\Bbb Z}
\define\r{\Bbb R}

\define\N{\Bbb N}
\define\jj{\Bbb J}

\define\Q{\Bbb Q}
\define\A{\Cal A}
\define\h{\Cal D}
\define\E{\Cal E}
\define\m{\Cal M}
\define\T{\Cal T}
\define\f{\Cal S}

\define\om{\omega}

\define\e{\varepsilon}
\define\va{\varphi }
\define\CB#1{\Cal C_b(#1)}
\define\st{\subset }
\define\al{\alpha  }
\topmatter
  \title
 Recurrent solutions of neutral  differential-difference systems
 \endtitle
\endkeywords
  \author
  Bolis Basit and Hans G\"unzler
\endauthor
 \abstract
{Results of Bohr-Neugebauer type are obtained for recurrent
functions : If $y$ is a bounded uniformly continuous solution of a
linear neutral difference-differential system with recurrent
right-hand side, then $y$ is recurrent if  $c_0 \not \subset X$ ;
also analogues and extensions to half lines are given. For this,
various subclasses ``$rec$" are introduced which are linear (the
set REC of all recurrent functions is not), invariant, closed etc.
Also, analogues of the Bohl-Bohr-Amerio-Kadets and Esclangon-
Landau results for REC are obtained.}
  \endabstract
 \endtopmatter
\rightheadtext{ recurrent solutions} \leftheadtext{Bolis Basit and
Hans G\"unzler}
 \TagsOnRight
\document
\pageno=1
 \baselineskip=16pt


 \head{\S  0    Introduction}\endhead

The concept "recurrent motion" has been introduced by Birkhoff
[14, p. 305, 311], generalizing the periodic case, to describe the
asymptotic behaviour of solutions of autonomous systems  of
differential equations (see also [16, vol. 1, survey, p. 282] ).
Since then recurrent has been extended to general dynamical
systems  [26], [34, p. 373-379], [2], [24], [25], [13, p.
195-199], with various meanings and applications in various
fields, even in number theory (see the references in [24, p.
117]).

  Though introduced before Bohr's almost periodic
functions, recurrent is more general than almost periodic, in fact
one has, for  $f :$  reals $  \r\to $ Banach space $X$, with all "
$\st$ " strict: (see (2.6), (2.9))

(0.1) \qquad $\{$periodic$\,\, f\} \st  \{$almost periodic$\,\,
f\} \st \{$almost automorphic$\,\, f\}$

   \qquad  \qquad  \qquad  $\st
          \{$ Levitan almost periodic$\,\, f\}  \st \{$recurrent $\,\, f\}$,

\noindent  using here and in the following "recurrent" in the
sense of [32, p. 80, Definition 2].

  Beginning with Bohr and Neugebauer [19], for the non-autonomous
case results of the type ``if $y$ is a bounded (uniformly
continuous) solution of a linear difference-differential system
with  right hand side  $f$ almost periodic, almost automorphic or
Levitan almost periodic, the $y$ is of the same type''  have been
obtained ([19], [17], [22], [6], [7]).

   Such results seem to be unknown for the recurrent
case, the main difficulty being that the sum of two recurrent
functions is in general not recurrent. Here we remedy this by
introducing in \S 3 certain subsets $rec \,\,V$  of recurrent
functions which are still linear, but large enough to have all the
other properties needed, so that the results of [7] (in a slightly
generalized form) can be applied:

  Bounded uniformly continuous solutions of even neutral
difference-differential systems are indeed recurrent, if the right
hand side is (Theorem 5.1).

  We also extend results of [8] to recurrent solutions on half
lines (Proposition 5.10/Corollary 5.11), and show that (weakly)
recurrent solutions have recurrent derivatives up to the order of
the system considered (Proposition 5.6).

 In \S 2, \S 4 we discuss explicit examples, e.g. of  recurrent $g_j$
with $g_1 + g_2$  not recurrent,  which show that the inclusions
between the various spaces considered are strict respectively that
our results are sharp.

  In \S 4 the analogue  of the
Bohl-Bohr-Amerio-Kadets result on the integration of almost
periodic functions [1, p. 55], [30] for recurrent functions [4] is
extended and obtained also for various subclasses of recurrent
functions as needed in \S 5.  Furthermore with Theorem 4.5 for the
first time a difference version of such results for recurrent
functions is obtained. Also, it is shown that though Theorem 4.5
is valid for $REC(\r_d,X)\cap C_b(\r,X)$, $\r_d= \r$ in the
discrete topology (Example 4.12 (a)), not even the classical
Bohl-Bohr analogue holds, that is already
 for $X = \cc$  the $\A = REC(\r_d,\cc)\cap C_b(\r,\cc)$
        does not satisfy
                $(P_b)$ (Definition 4.1, Example 4.12 (b)).

 \head{\bf\S 1. Notation and  Definitions }\endhead

In the following $\jj$ will always  be an interval  of the form $
[\al, \infty)$ for some  $\alpha \in \r$, $\r_+ =[0, \infty)$,
$\r^+=(0, \infty)$, $\N= \{1, 2,\cdots \}$ and $\N_0=\{0 \} \cup
\N$.
 Denote by $X$ a real or  complex Banach space, with scalar field  $\r$ or
 $\cc$. We write $c_0 \not \st X$ if $X$ does not contain a subspace isomorphic to the Banach space $c_0$.
 If $f$ is a  $ X$-valued function defined on  $T\in \{\r, \jj,$ abelian group $G \}$,
then   $f_s$, $\Delta_s f$ will stand for the functions defined on
$T$  by $f_{s}(t) = f(t+s)$,
 $\Delta_s f(t) = f_s (t) - f (t)$ for all $s\in T $ with $s+T \st
 T$,
  $|f|$ will denote the function $|f|(t):= ||f(t)||$  for all $t\in T$ and $||f|| _{\infty} :=
  \text {  sup}_{x\in T} ||f(x)||$.
 If $f \in L_{loc}^1 (T, X)$ with $T= \r$ or $\jj$, then $Pf$ will denote the
indefinite integral defined by $Pf(t) = \int_{\alpha}^{t}
f(s)\,ds$ (where   $\al =0$   if $T=\r$, all integrals are
Lebesgue-Bochner integrals (see [3, pp. 6-15], [23, p. 318],
   [29, p. 76]), similarly for measurable).

     In $X^{T}$,  the spaces of all constants, continuous, bounded continuous, uniformly continuous,
     continuous with relatively compact range and continuous with relatively weakly compact range
  will be
denoted by

$X$, $C(T,X)$, $C_b(T,X)$, $C_u(T,X)$, $C_{rc}(T,X)$,
$C_{rwc}(T,X)$;

  $C_{ub}(T,X)= C_u(T,X)$
 $\cap\, C_b(T,X)$, $C_{urc}(T,X)=C_u(T,X)\cap\, C_{rc}(T,X)$
etc.

  For $C_{loc} (T,X)$ see (2.3), for $C_{b, wscp} (T,X)$ see
(4.2).

\noindent If $\A \st X^T$, $\A_b$, $\A_u$, $\A_{rc}$, $\A_{rwc}$,
$\A_{ub}$, $\A_{urc}$ mean $\A \cap C_b$, $\A \cap C_u$ etc.
\smallskip

(1.1) \qquad  $\m\A : = \{f \in L^1_{loc} : M_h f (\cdot) : =
(1/h)\int^h_0 f(\cdot+s)ds \in \A$,\,\, all\,\, $h > 0 \}$,

\qquad\qquad  \, $\m ^0 \A : = \A \cap L^1_{loc} (T,X)$,
 \,\,\, $\m^k \A : = \m(\m^{k-1}\A)$,\,\,\,\,\,\, $T=\r$ or $\jj$.

\smallskip

 \noindent Usually  $ \A  \st   \m\A  \st \m^2\A  \st \cdots $ even strictly ([9,
section 2], Example 4.12 (c)).

\noindent  For  $\A  \st  X^T$  in the following "$\A$
positive-invariant" etc means:

 $Positive$-$invariant$: translate  $f_a  \in   \A$  if  $f  \in
 \A$
and  $0\le a <\infty$.

 $Invariant$ :  $f_a  \in   \A$ if  $f  \in  \A$, $a  \in  T$,
$T= \r$ or $G$.

 $Uniformly\,\, closed$ : $f_n  \in  \A$, $n \in  \N$, and $f_n \to
 f$
uniformly on $T$ implies $f \in \A$.

$(\Delta)$ :  $f  \in   L^1_{loc}(T,X)$, $\Delta_h f  \in   \A$
for \,\,\, $h>0$ \,\,\,\, implies \,\,\,\, $f - M_k f \in\A$
                                 \,\, for \,\, $k > 0$,

\qquad\qquad     $ T = \jj$ or $\r$;\,\,\, $(\Delta)$ has been
found useful in [8],
     [9], [10].

\noindent Further definitions and function spaces:

 $  (\Gamma)$: Definition 3.1;\,\,\,
$(L_U)$, $(P_U)$: Definition 4.1;\,\,\, $(P'_U)$: after (5.6).

    $\lambda$-class: Definition 3.1.

 $AP\,\, = \{$ almost
periodic $f\}$, \qquad $ (B) AA = \{$almost automorphic  $f\}$,

 $LAP =
    \{$ Levitan-almost periodic $f\}$: (2.5) and after.

$ REC = \{$ recurrent  $f\}$,\,\,\, $MIN =\{$ minimal $f\}$:
Definitions 2.1, 2.2, (5.6).

 $rec \, V$, $rec_b \, V$ $rec_u \, V$, $\cdots$: Definition 3.2,  before (3.1); $O(f)$: after (2.4).

$ \E = \{$ ergodic $f\}$: after Proposition 4.2;\,\,\, $\T\E$:
(5.3);\,\,\, $O_q$ : (5.7).

$\r_d$, $G_d$: $\r$ respectively $G$ in the discrete topology;
$O_d (f)= O(f| G_d)$.

\head{\bf \S 2.   Recurrent Functions}\endhead

  As mentioned before, "recurrent" can have various meanings; here we
use it in the original sense of Birkhoff [14] or Levitan and
Zhikov [32, p. 80 Definition 2], see also Flor [24, Definition 1],
without the "uniform" :

\noindent {\bf {Definition 2.1}}. If $G$ is a topological abelian
group and $X$ a Banach space,
     an  $f : G \to X$  will be called $recurrent$ if it is continuous and if for
     each compact $K\st G$ and $\e > 0$  the set of $\e$-$K$-periods
     $E(f,\e,K)$ is relatively dense;
      here

\smallskip

(2.1) \qquad  $E(f,\e,K) : = \{\tau\in G : ||f(t+\tau)-f(t)|| \le
\e$ for $t\in K\}$,

\smallskip

  \noindent   an  $E\st G$   is called $relatively\,\, dense$ if there exists a compact $L\st
     G$
     with  $E + L = G$.

\smallskip

(2.2) \qquad  $ REC(G,X) : =  \{ f \in X^G  :  f $ recurrent $\}
$.

\smallskip
If $f \in C(G,X)$, then for each $\e >0$ and each compact $K\st
G$, $G$ locally compact, there is a compact neighbourhood  $W\st
G$ of $0$ such that

(2.3)\qquad $E(f,\e,K+W)+W \st  E(f,2\e,K)$.

\noindent So $REC (\r,X)=RC (\r,X)$ of [10, before (1.6)]:   $For
\,locally \, compact \, abelian \, G, \, an $

\noindent $f \in C(G,X)$ $\,\, is \,\, recurrent \,\, if \,\, and
\,\,
     \,\, only \,\, if\,\, for\,\, each\,\, \e > 0\,\,\, and\,\, compact\,\, K  \st G \,\,$ $ there \,\,exists
     \,\,a$
 $
       finite \,\, F_{\e,K}  \st  G \,\, with \,\, E(f,\e,K) + F_{\e,K} = G $; so

       \qquad\qquad\qquad $REC (G,X)\st REC (G_d,X)$,

       \noindent with $G_d = G$ in the discrete topology (see
       (2.13)).

$REC(G,X)$ is closed with respect to uniform convergence on G and
invariant, i.e. translate $f_s \in REC$  if  $f \in REC$ and $s
\in G $.

 If  $G =\r$, then  $K = [-n,n]$ suffice, $E\st \r$  is relatively
dense if and only if there is $L \in [0,\infty)$ such that $E \cap
[a,a+L] \not=  \emptyset$  for each $a \in \r$.

  If one equips $C(G,X)$ with the topology of locally uniform
convergence (i.e. on each compact set), one gets a locally convex
linear space which we denote by

\smallskip

(2.4) \qquad      $Y  =  C_{loc}(G,X ) $; $C_{loc}(G,X)$ is
complete if $G$ is locally compact.

\smallskip

 \noindent For  $f \in C_{loc}(G,X)$, the  $orbit\,\,$ $O(f) : =$ closure in $Y$ of
$\{f_s : s \in G\}$;  $O(f)$ is invariant.

\noindent{\bf {Definition 2.2}}. An $f \in C_{loc}(G,X)$  is
called $minimal$, if its orbit  O(f) is
    minimal, i.e.    $\emptyset \not = U \st O(f)$ with $U$ invariant and closed in $C_{loc}
    (G,X)$  always implies $U = O(f)$;

\smallskip

(2.5) \qquad    $MIN (G,X) : = \{f\in C_{loc}(G,X) : f$ minimal
$\} $.

\smallskip

For general abelian topological $G$ and Banach space $X$ one
always has, with all " $\st$ " being strict if $G =\r$,

\smallskip

(2.6)\qquad  $AP\st  BAA\st AA\st LAP\st REC\st MIN $.

\smallskip

\noindent Here AP, BAA, AA and LAP denotes respectively the linear
space of Bohr-Bochner almost periodic (= ap), Bochner almost
automorphic , Veech almost automorphic continuous (= aa) and
Levitan ap functions [39, p. 18], [18], [37], [35, Definition 2]
and [10, \S 3] ($RC$  of [10] is our $REC$ by (2.3)); for the
first 3 strict " $\st$ " see [10,(3.3)] valid also for G as here,
the next " $\st $ " follows from [35, Satz 1, Satz 2], the last
from [24, Satz 5, "f(G) totally bounded" is not needed]; $REC$
strictly $\st MIN $ is shown in [24, first example p. 127], $LAP$
strictly $\st REC$ by Example 2.10 below (see also (2.9)).

  If $G$ is $\sigma$-compact and locally compact, then $C_{loc}(G,X)$ is
metrizable; for $G =\r$ one can use e.g. ([2, p. 145]

(2.7) \qquad $d(f,g) : =$ sup $_{n \in \N} $ min $\{1/n,\,\,$ sup
$_{|t|\le n} ||f(t)-g(t)||\} $;

\noindent then a continuous $f : G \to X $ is recurrent if and
only if for each $\e > 0 $

\smallskip

(2.8)  \qquad $\{ \tau\in G :  d(f, f_{\tau}) \le \e \}$ is
relatively dense.

\proclaim {Theorem 2.3 (Flor)} If G is locally compact and  $f \in
Y : = C_{loc}(G,X)$,
     the following five conditions are equivalent :

(a)    $O(f)$ is compact and minimal in Y,

(b)    $O(f)$ is minimal in $Y$ and  $f \in C_{urc}(G,X)$,

(c)    $f$ is recurrent and  $\{f_s : s \in G \}$ totally bounded
in $Y$,

(d)    $f$ is recurrent, uniformly continuous, with  $f(G)$
relatively compact in $X$.

(e) $f$ is uniformly continuous and to each net $(r_i)$ from $G$
there is a subnet

\qquad $\sigma =
      (s_j)$  and a net $\tau = (t_k)$ such that  $\sigma f : = lim f_{s_j}$ exists
      locally uniformly

   \qquad   and $\tau(\sigma f) = f$ locally uniformly.

\endproclaim
 With (2.6) and [10, (3.3), Proposition 3.1] this implies, for locally compact $G$,
any $X$,

\smallskip

(2.9)  \qquad $AP\st BAA_u = AA_u = LAP_{urc} \st REC_{urc} =
MIN_{urc}$

  \qquad\qquad  \qquad\qquad                         $ \st REC_{ub}\st REC_b\st REC $,

\smallskip

again all " $\st$ " are strict by  Example 2.7, Example 2.10,
Remark 2.5 (i), [10, Example 3.3], [32, p. 58, 4.].

\proclaim {Corollary 2.4} If  $f \in REC_{urc} (G,X)$ with $G$
locally compact, then

\qquad \qquad $O(f)
                 \st REC_{urc}(G,X)$.
\endproclaim

\proclaim {Remarks 2.5} (i). If $G$ is arcwise connected, then

(2.10)\qquad $REC_u(G,X) =REC_{ub}
       (G,X)$

 \noindent arguing as in [10, Proposition 3.1] for the case $G =\r$. However,
      $REC_{ub}
       (\r,X) \not\st C_{rc}(\r,X)$  with
       $\phi$ of [5, Lemma 5.1], $x_k =$
        basis $e_k$ in  $l^p$
       in (5.1)
there,  $\pi(\phi)  \in  LAP_{ub}(\r,l^p)$, $1\le p< \infty$,
$\not \in C_{rc}$, with $\pi$ of Example\, 2.7 below).

   (ii)  Theorem 2.3 becomes false if one drops either "compact" or              "
minimal"  in (a), or in (b) ([24, examples p. 127]).
\endproclaim

\demo{Proof of Theorem 2.3} (See also [13, p. 198-199, 8.8
Theorem], [34, p.373-379] for G = continuous 1-parameter
         group.)

  $a \Rightarrow b$ : Flor [24, Satz 8,
Satz 3, Satz 4] .

 $b \Rightarrow d$: Flor [24, Satz 8].

 $d
\Rightarrow c$: Follows from a general Arzela-Ascoli result [24,
proof
     of Satz 2], [23, p. 266/267] :

\smallskip

(2.11) \qquad $ f  \in  C_{urc}(G,X)$ implies  $O(f)$ compact in
$C_{loc}(G,X)$,

 \qquad  \qquad  \qquad  \qquad  and
                                         $ O(f)  \st  C_{urc}(G,X)   $.

\smallskip

$c \Rightarrow a$ :  [24, Satz 5], totally bounded = relatively
compact in
      the complete Y [31, p. 36].

$a \Leftrightarrow e$: [25, Lemma 1, Definition 2].
      $\square$ \enddemo

\demo {Proof of Corollary 2.4} : $ d \Leftrightarrow a $ and $O(g)
= O(f)$ if $g \in O(f)$ minimal . $\square$
\enddemo

In section 3 we need

\proclaim{Theorem 2.6 (Flor)} If $G$ is an abelian locally compact
group, $ f, g  \in
       C_{urc}(G,X)$, $f$ almost automorphic  $\in  AA(G,X)$, $g$ recurrent  $\in
        REC(G,X)$,  then  $(f,g)$ is also recurrent,  $\in  REC(G,X^2)$.
\endproclaim

So in the terminology of [32, p. 10]  or [25, Definition 6], aa
$f$ are "absolutely recurrent" (but restricted to $C_{urc})$; for
ap $f$ this can be found also in [32, p. 106].

\demo{Proof} (2.11) and a contradiction argument yield, for
$\va\in C_{urc}(G,X)$, with orbit $O(\va)$  in $C_{loc}(G,X)$,

\smallskip

(2.12) \qquad If $(g_j)$ is a net from $O(\va)$, $h \in  X^G,
g_j\to h$ pointwise on $G$,
                  then

                  \qquad\qquad \qquad $g_j\to h$  locally uniformly in $G $.

\smallskip

\noindent So if  $G_d : =  G$ with the discrete topology  and $U
\st O(\va)$, $U$ is closed in $C_{loc}(G,X)$ if and only if  $U$
is closed in $C_{loc} (G_d,X)$, i.e. with respect to pointwise
convergence. With $O : = O(\va)$,  $O_d : = O(\va|G_d)  = $ orbit
of $\va$ in  $C_{loc}(G_d,X)$, obviously  $O  \st  O_d $; (2.12)
gives $O_d \st O$, so $O = O_d $. This implies, for any  $\va \in
C_{urc}(G,X)$

\smallskip

 (2.13)\qquad  $O$ is minimal in $C_{loc}(G,X)$ if and only if  $O_d$  is
minimal in $C_{loc}(G_d,X) $.

\smallskip

\noindent $g \in REC_{urc}(G,X)$, Theorem 2.3 and (2.13) yield
therefore $g \in MIN_{urc} (G_d,X)$.

  Now by Satz 4 of  [24], in our notation,

   $f \in
AA(G_d,X)$ and  $g \in  MIN_{rc}(G_d,X^2)$  imply $(f,g)\in
MIN_{rc}(G_d,X^2)$.

  Since $(f,g) \in  C_{urc}(G,X^2)$, (2.11)), $O = O_d$ of above, (2.13)
and Theorem 2.3 for $X^2$  give $(f,g) \in REC(G,X^2)$. $\square$
\enddemo

\proclaim{Example 2.7}  Explicit almost automorphic uniformly
continuous (then
    Bochner aa  =  Veech aa)  $f : \r\to \cc$   which is  not Bohr almost
    periodic:

     For  $n \in \z$  define  $\va(n) : = (1 + e^{in}) /
    |1 + e^{in}| $;  then  $\va \in  AA(\z,\cc)$  but  $\va  \not \in
    AP(\z,\cc)$
    [13, p. 192, Corollary 7.15, Example 7.16a (without (1.9)].
\endproclaim

  If  $\pi(\va)$  is the continuous piecewise linear  extension of
$\va : \z\to \cc$
    to $\r$ of [5, Th\'{e}or\`{e}me 5.3], then  $f : =  \pi(\va) \in  BAA_u
     (\r,\cc)  =  AA_u(\r,\cc) , =  LAP_{urc}(\r,\cc)$ of [10, (3.3)], using
     [10, Theorem 3.2(e)], valid also for $\z$ . With Bochner's characterization or ap,
      this f  is not in $AP(\r,\cc)$.

   A more involved first such example is due to Veech [37,
Theorem
     6.1.1 and Remark].

\proclaim{Example 2.8} $REC_{urc}(\r,\cc)$ and so $REC(\r,\cc)$
are not linear spaces.

\endproclaim

\demo{Proof}   With Kronecker's approximation theorem [28, p. 436,
(d)]
   one can find $\psi_1, \psi_2 \in  O(\va)$ with $\va$ of Example 2.7  with

  $ \psi_1(n) = \psi_2(n) =  (1 - e^{in})/|1 - e^{in}|$  if  $0\not =n
\in \z$,
      $\psi_1(0) = i$, $\psi_2(0) = -i $;

 \noindent  since  $\va \in  REC_b(\z,\cc) = REC_{urc}(\z,\cc)$ by (2.6) with
$G = \z$ in the
      discrete topology, Corollary 2.4 yields $ \psi_1, \psi_2 \in REC_b (\z,X)$.
      With $g_j : = \pi(\psi_j)$  as in Example 2.7, $ g_j \in  REC_{urc}(\r,\cc)$, using the
      uniform continuity of the $\pi(\psi)$ for $\psi$  bounded.
    Obviously  $g_1 - g_2$  is not recurrent showing $REC_{urc}(\r,\cc)$ is not linear.
    (We note also that
       $g_1 + g_2 \not   \in  REC(\r,\cc)$.  Since  we do not use
       this, we omit the proof.) $\square$
       \enddemo
    A first such example seems to be due to Auslander and Hahn
[2,  Example 2.1].

\proclaim{Example 2.9} The function $g= e^{i\,t^2}$ belongs to
$MIN_b (\r_d,\cc)=REC_b(\r_d,\cc)$ and  $MIN_b (\r,\cc) $, but not
to $REC_b(\r,\cc)$ so ($\r_d$: $\r $ in the discrete topology)

(2.14)\qquad  $REC_{rc}(\r,\cc)$ strictly $\st MIN_{rc}(\r,\cc)$,

       \qquad \qquad \,\,\,\,            $REC_{rc}(\r,\cc)$ strictly $\st REC_{rc}(\r_d,\cc)  \cap  C(\r,\cc) $.

     \noindent  This $g$ also shows  $MIN_{rc} \not \st C_u$;  $MIN_u \not \st C_b$ follows
       from [24 , example p. 127].
\endproclaim
\demo{Proof} $f\in O_d(g):= O(g| \r_d)$ implies $f=c\cdot g\cdot
\chi$, where $|c|=1$ and $\chi \in \widehat{\r_d}$. By  [13, 6.12
Proposition, p. 53] and Theorem 2.6, $g\in MIN (\r_d,\cc)\cap C_b
(\r,\cc)$ (see also [13, Exercise 8.18, p. 202]).  $g\not \in REC
(\r,\cc)$ with Definition 2.1. $\square$
\enddemo

\proclaim{Example 2.10} There are   $g \in REC_{urc}(\r,\cc)$
which are not Levitan ap.
\endproclaim
        \demo{Proof}  One of $ g_1, g_2$ of Example 2.8  cannot be in
       $LAP(\r,\cc)$,  for then, $LAP$ being linear, also  $g_1 - g_2 \in  LAP\st   REC$, contradicting Example
       2.8. $\square$
\enddemo

\proclaim{Example 2.11}  $\A$  strictly  $\st  \m\A$  for $\A =
REC(\r,X)$, $REC_b$, $REC_u$,
     $REC_{urc}$, $rec_{xy} V $ of Definition 3.2 , there is $f  \in  C_b(\r,\r)$  which is not
    recurrent (not even Poisson stable [32, p. 80 Definition 1], also not Stepanoff
     $S^1$-almost periodic), but  $Pf$  is almost periodic  and so  $f  \in \m AP
     \st  \m REC_{urc}$:

      $f =  \sum_{n=2}^{\infty} h_n$,\qquad $h_n$ periodic  with period $2^{n+1}$,

     $ h_n (t)=  0  $, $\,\,\, t \in [-2^n, 2^n  - 1]$,\qquad
       $h_n(t)= sin(2^n \pi (t - 2^n))$, $\,\,\, t\in [2^n  -1,2^n]$;
     $V = \{M_h f : h>0 \}$ for $ rec_{xy} V $.
\endproclaim

\proclaim{Proposition 2.12} $REC(\r,X)\st \m REC(\r,X)$.
\endproclaim
\demo{Proof} For   $g \in L^1_{loc}(\r,X)$; $0 \not= h \in \r$, $n
\in \N$, with
      $E(g,n)  : =     E(g,1/n,[-n,n])$   of  (2.1)
one has, by integrating the inequality in (2.1), with $[h] : =$
largest integer $\le h$,

(2.15)   \qquad $ E(g, n + [ |h| ] + 1)    \st E(\Delta_h Pg, n)$.
\qquad\qquad  $\square$
\enddemo

\head{\S 3  Construction of a recurrent $\lambda$-class}\endhead

  To apply the results of [7] to recurrent solutions of differential-difference systems,
   a $\lambda$-class (Definition 3.1 below) of such
functions is needed. $REC(\r,X)$ or $REC_{urc}(\r,\cc)$ cannot be
used, since these spaces are not even linear (Example 2.8), also
the $r(\va,\r,X)$ of [6] is not good enough. Here we remedy this.

\proclaim {Definition 3.1} An  $\A  \st  X^G$  is called a
$\lambda$-class ([6, Definition 1.3.1], [7, p. 117]) if $\A$ is an
invariant closed linear subspace of $C_{ub}(G,X)$ containing the
constants which satisfies

$(\Gamma)$ \,\,   $f  \in  \A$,  $\omega  \in  \widehat {G}$ imply
$\omega f\in  \A$.
\endproclaim
  \noindent  Here $G$ is an abelian topological group, $X$ now a complex Banach
    space, $C_{ub}$ a Banach space with  $||f||_{\infty} : =  sup_{t\in G} ||f(t)||$,  $\A$ invariant
    means translate  $f_s  \in  \A$  if $f  \in  \A$, $s  \in  G$, and $\widehat {G}$  contains all
    continuous bounded characters  $\omega : G\to\cc$,  $\om (s+t)=\om(s)
    \om(t)$ for $s, t  \in  G$ ($\widehat {\r} = \{ \gamma_{\omega} : \om  \in\r\}$ with
    $\gamma_{\omega}(t) = e^{i \om\, t}, t  \in\r $ [28, p. 367, (e)]).

 \noindent To get such a $\lambda$-class $\st REC$, for  fixed $V$ with
  $\emptyset \not = V\st  REC_{xy}(G,X)$,

  $ REC_{xy} \in \{REC, REC_{b}, REC_{u}=REC_{ub},
     REC_{urc}, \cdots\}$ (see Remark 2.5 (i)),

  \noindent  we introduce

\smallskip

(3.1) \qquad $\frak{V}: = \{U \st REC_{xy}(G,X) : V\st U$, $U$
satisfies (3.2)$\}$, \,\, with

\smallskip

(3.2) \qquad $n \in \N$, $f_1,...,f_n \in U$ implies
$(f_1,...,f_n)\in REC(G,X^n) $.

\smallskip

\noindent Then $\frak{V}\not =\emptyset$ if ${V}$ satisfies (3.2),
$\frak{V}$ is partially ordered by "$\st$" such that any chain
from $\frak{V}$ has  the union over the chain as a sup $\in
\frak{V}$. So Zorn's lemma can be applied ([31, p. 8]), there
exists a maximal element in $\frak{V}$.

\proclaim{Definition 3.2} For fixed $V$ with $\emptyset \not =
V\st REC_{xy}(G,X)$ and satisfying (3.2), we denote by $rec_{xy}\,
V$ one
    such maximal element of  $\frak{V}$,

  \noindent   $REC_{xy}\in \{REC, REC_{b}, REC_{u}=REC_{ub}, REC_{rwc}, REC_{rc},REC_{urwc}, REC_{urc}\}$.
\endproclaim
\proclaim{Proposition 3.3} If  $\emptyset \not = V\st
REC_{xy}(G,X)$ with (3.2) and $G$, $X$, $xy$ are as in Definitions
3.1/3.2,
        $rec_{xy}\, V$ of Definition 3.2 is  linear, invariant, uniformly closed, contains
                        the constants, with  $V \st rec _{xy}\, V \st REC_{xy}
        (G,X) $.
\endproclaim

\demo{Proof}  (a)  $L : = rec_{xy} V$ is linear:  If  $g,h,$
$f_1,...,f_{n-1}\in L$, then
 $(f_1,...,f_{n-1},g,h) \in REC$; this implies $g+h \in REC_{xy}$ and
$(f_1,...,f_{n-1},g+h)$, $(f_1,...,f_k,g+h,...,g +h) \in REC$  if
$0\le k< n-1$, also in arbitrary order. Since $L$ is maximal, $g+h
\in L$. $c g \in L$  for complex $c$ follows similarly.

 (b) $L$ is uniformly
closed: If $g_k \in L$  with  $||g_k - g||_{\infty}\to 0$ and
$f_1,...f_{n-1}\in L$, then $(f_1,...,f_{n-1},g_k) \in REC $ for
each $k$ implies $(f_1,...,f_{n-1},g) \in REC$  since $REC(G,X^n)$
is closed with respect to  $|| \cdot ||_{\infty}$. As in  (a) ,
this gives $g \in L $.

(c) L is invariant :  If  $f_1,...,f_{n-1},g \in L$, $K$ compact
$\st G$, $\e
> 0 $, $a \in G$  there is $\delta > 0$ with  $E : =
E((f_1,...,f_{n-1},g), \delta, K \cup (a+K)) \st E(g_a , \delta,
 K)$ and $\st E(((f_1,...,f_{n-1}, \e, g_a),K)$. Since $E$ is
relatively dense, $(f_1,...,f_{n-1},g_a) \in REC$ and then $g_a
\in L$. $\square$
\enddemo

\proclaim {Corollary 3.4} If  $\emptyset \not = V\st
REC_{urc}(G,X)$ with (3.2), then $rec_{urc} \,V$ is a
$\lambda$-class.
\endproclaim
\demo{Proof} $L$ of the proof of Proposition 3.3 satisfies
$(\Gamma)$ : If $f_1,...,f_{n-1},g \in L$, $\omega \in \widehat
{G}$, $0\not = a \in X^n$, $K$ compact $\st G$, $\e
> 0$, since $f$ and $\omega$ are bounded, there is $\delta > 0$ with

  $E :
= E((\omega\, a,(f_1,...,f_{n-1},g)), \delta, K)\st
E((f_1,...,f_{n-1},\omega g), \e, K)$.

\noindent Since $E$ is relatively dense by Theorem 2.6, one gets
$(f_1,...,f_{n-1},\omega g) \in REC$ and thus $\omega g \in L$.
$\square$
\enddemo
\proclaim {Remark 3.5} The $rec_{urc} \, V$  of Corollary 3.4
contains also all uniformly continuous almost automorphic
functions,

 $AP(G,X) \st AA_u(G,X) = LAP_{urc}(G,X)\st rec_{urc} \,
V$.
\endproclaim
\noindent (Proof as for Corollary 3.4 above.)

 \proclaim {Corollary 3.6} If  $\emptyset
\not = V\st REC_{xy}(\r,X)$ with (3.2), $xy$ as in Definition 3.2,
then $\A =rec_{xy} V$ satisfies $\A \st \m \A$;  $\A$ satisfies
$(\Delta)$ if $\A \in \{ rec_{u}\,V, rec_{urc}\, V\}$.
\endproclaim
\demo{Proof} $\A  \st  \m\A$ follows with (2.15), the cases
    $rec_{rc}$, $rec_{urc}$ with [10, Proposition 2.2 (i)].
If  $f  \in  C_{rwc}$, then $M_h f  \in  C_{rwc}$  by
Krein-Smulian [23,
        p. 434, V.6.4 Theorem] as in the proof of Proposition 2.2(i) of [10] .

     $(\Delta)$ follows by Proposition 3.3 and [10,
Theorem 2.4]. $\square$
\enddemo
\proclaim{Proposition 3.7}  Assume  $f\in\m REC _{xy}(\r,X)$, $xy=
ub$ or $urc$. Then

(i) $f\in \m rec_{xy} \, V$, where $V := \{M_{h}f: h>0,  h\in
\Q\}$ satisfies (3.2).

(ii) $f\in \m rec_{xy} \, V$, where $V := \{M_{1/2^n}f:  n\in N\}$
satisfies (3.2).
\endproclaim
\demo {Proof}(i) Let  $h=p/q$, $k= r/q$ with $p, q, r \in \N$.
Then $h M_h f= q (M_{1/q} f + (M_{1/q} f)_{1/q}+\cdots+ (M_{1/q}
f)_{(p-1)/q}) \in rec_{xy} \{M_{1/q} f\}$ by Proposition 3.3.
Similarly $k M_k f = q (M_{1/q} f + (M_{1/q} f)_{1/q}+\cdots+
(M_{1/q} f)_{(r-1)/q}) \in rec_{xy} (M_{1/q} f)$. Since  $
U=rec_{xy} \{M_{1/q} f)\}$ satisfies (3.2), one gets $(M_h f, M_k
f)\in REC_{xy} (\r,X^2) $. In the same way, one gets $(M_{h_1} f,
\cdots,  M_{h_n} f)\in REC_{xy} (\r,X^n) $ for all $h_1, \cdots,
h_n \in \Q^+$. This proves that $V := \{M_{h}f: h>0,  h\in \Q\}$
satisfies (3.2).

 Since $f\in \m REC _{ub}(\r,X)$, by [7, Theorem 3.3] one gets $Pf \in C_u
 (\r,X)$. It follows  that if $k_n \to k$  as $n\to \infty$, then,
 since  $M_k f \in C_{ub} (\r,X)$,

 (3.3) \qquad $M_{k_n} f  \to M_k f $ uniformly on
 $\r$.

 \noindent  Since $\Q$ is dense in $\r$ and $\A= rec_{xy} \, V$ is  uniformly closed, by (3.3) and the above, one gets
 $M_k f \in \A$, $k>0$.

(ii) Similarly as in (i). $\square$
\enddemo

\proclaim{ Corollary 3.8} $\m REC_{urc}(\r,X)$  satisfies
$(\Gamma)$.
\endproclaim
\demo {Proof}  If   $f  \in \m REC_{urc}(\r,X)$, with $f \in  \m
rec_{urc}\,  V$  of Proposition 3.7, $rec_{urc}\, V $ is a
$\lambda$-class with $(\Delta)$ by Corollaries 3.4/3.6, then $\m
rec_{urc}\, V $  has $(\Gamma)$ by [12, Proposition 5.1(ii)],
yielding $(\Gamma)$ for $\m REC_{urc} $ since
 $\m
rec_{urc}\, V\st \m REC_{urc} $. $\square$
\enddemo

\proclaim{Proposition 3.9}  If  $\A=REC _{xy}(\r,X)$ with $xy$ as
in Definition 3.2, then

(i)\,\, $\A\st \m \A$.

(ii)\,\, $\A$ satisfies $(\Delta)$ if $xy\in \{u, urc, urwc\}$
(see also Corollary 4.6).

(iii) If  $f\in \m \A$ with $xy=u$ or $urc$, $n\in \N$, $h_j >0$,
then

\qquad $(M_{h_1} f, \cdots, M_{h_n}) \in REC _{xy}(\r,X^n)$.
\endproclaim

\demo{Proof} (i)  as for Corollary 3.6.

(ii) Let $f\in L^1_{loc} (\r,X)$ and $\Delta_h f \in REC
_{urc}(\r,X)$ for $h>0$. By  [7, Theorem 3.3], $f\in C_u (\r,X)$
and by [9, Proposition 3.1], $f-M_k f\in  C _{urc}(\r,X)$. This
implies

(3.4)\qquad $f -M_k f = \lim_{n\to \infty} \sum_{j=1}^n (1/n)
(f-f_{(kj/n)})$ uniformly on $\r$.

\noindent As $f-f_{(kj/n)}= f-f_{k/n}+ (f-f_{k/n})_{k/n}+\cdot +
(f-f_{k/n})_{(j-1)k/n}\in rec_{urc} \{f-f_{k/n}\}$ for $j=1,
\cdots, n$, $ \sum_{j=1}^n (1/n) (f-f_{(kj/n)})\in rec_{urc}
\{f-f_{k/n}\} \st REC _{urc}(\r,X)$ for each $n\in \N$ and $k>0$
by Proposition 3.3. By (3.4) and the remark after (2.3), $f-M_k
f\in REC_{urc} (\r,X)$ proving $(\Delta)$.

Case $xy=u$ respectively $urwc$: Similar, with Remark 2.5 (i)
respectively [7, p. 119 below] or Lemma 4.6 below.

(iii) follows with Proposition  3.7 and Definition 3.2. $\square$
\enddemo

\head{\S 4  Analogues to the Bohl-Bohr-Amerio-Kadets
Theorem}\endhead

  The BBAK-theorem of the title says that the indefinite integral of an
almost periodic  $f : \r\to X$  is again ap if it is bounded and
$c_o$ is not
 isomorphic to a subspace of $X$, or   $c_0 \not \st X$ for short
 ([30], [32, p. 84],  [3, Theorem 4.6.11, p. 298], [39, p. 99]).

  For the study of (generalized) ap solutions of differential
equations, a generalization of this is needed : If  for  $F :
\r\to X$ all differences $\Delta_h F = F_h - F$ are ap, $F$ is
bounded and $c_0 \not \st X$, then $F$ is ap([33], [21], [27],
[4]).

 To treat these and similar cases in a unified way, the following
Loomis-condition [33] has been introduced ([8, p. 677, 688], [10,
\S 6, D]):

\proclaim{Definition 4.1} If  $\A, U  \st  X^{\r}$, we say that
$\A$ satisfies $(L_U)$ if
        $F  \in  U$ with all differences  $\Delta_h F \in \A$  for  $h>0$ imply $F\in\A$.
        $\A$  satisfies $(P_U)$ if  $f \in \A$  with indefinite integral  $Pf \in U$  implies  $Pf \in \A$.
        $(P_b)$ means $(P_{C_b})$, similarly for $(L_b)$.
\endproclaim

  The BBAK-theorem  says then that $AP(\r,X)$ satisfies $(P_b)$ if
$c_0 \not\st  X$.

  As a first result one has

\proclaim{Proposition 4.2} For  $\emptyset\not = V \st
REC_{xy}(\r,X)$ with (3.2) the $rec_{xy} \, V$ of Definition 3.2
 and  $REC_{xy} (\r,X)$ satisfy
     $(P_{\E})$  and  $(L_{\E})$, $\E = \E(\r,X) = \{ergodic\, f\}$, $xy =u$ or $urc $.
\endproclaim
\noindent  Here an  $f \in L^1_{loc}(\r,X)$  is called  $ergodic$
([7, p. 117]) if there is  $m \in X$  with  $\lim_{T \to \infty}$
sup $_{x \in \r} ||(1/2T) \int^T_{-T} f(t+x)\, dt -  m ||    = 0$.

\demo{ Proof}  By Proposition 3.3,  $rec_{xy} V$ is closed linear
invariant $\st C_{ub}(\r,X)$, so by [9, Proposition 3.1] the
$rec_{xy} V$ satisfies $(\Delta)$.  Since $rec_{xy}\, V$ also
contains the constants, [8, Corollary 4.2] yields $(L_{\E})$.
$(L_{\E})$ implies here $(P_{\E})$ since for $ \A= rec_{xy} V$ one
has $\A\st \m \A$ by Corollary 3.6.

 If  $f  \in  L^1_{loc}$ with $\Delta_h  \in REC_{xy}$,  $\st REC_u = REC_{ub}$ of
      (2.10), by [9, Proposition 1.4] one has  $f  \in  C_u(\r,X)$.

   Now $V : = \{\Delta_r f : r  \in  \Q^+\}$ satisfies (3.2) : If
$r_1,...,r_q  \in
       \Q^+$,  $r_j = m_j/n$ and so  $\Delta_{r_j} f  \in  rec_{xy} \{\Delta_{1/n} f
       \}$
       by the formula after (3.4) with Proposition 3.3; therefore

      (4.1)\qquad   $(\Delta_{r_1} f,...,\Delta_{r_q} f)   \in   REC_{xy}(\r, X^q) $.

       If  $h$ real  $>0$ choose  $(r_m)  \st  \Q^+$  with $r_m \to h$ ; then $\Delta_{r_m} f
      \to  \Delta_h f$   uniformly on $\r$,  so  $\Delta_h f  \in  rec_{xy} V $ by
      Proposition 3.3 . $(L_E)$ for  $rec_{xy}$,  $\st REC_{xy}$, gives $(L_E)$ for $REC_{xy}$,
      then $(P_E)$ with Proposition3.9(i).
$\square$ \enddemo

For a direct analogue to the BBAK-theorem
 we need first
   generalizations of Lemma 1 and Lemma 2 of [4]:

(4.2) \qquad $C_{b,wscp}(\r,X) : = \{f \in C_b(\r,X): $

 \qquad \qquad\qquad\qquad \qquad\qquad\,\, $f(\r)$ relatively weakly  sequentially  complete $\}$,

\noindent where  $M  \st  X$ $relatively \,\, weakly \,\,
sequentially \,\, \,\, complete$ means to any sequence $(x_n)_{n
\in \N}$ with  $x_n \in M$ which is weakly Cauchy there is $a \in
X$ with $x_n \to a$ weakly.

\proclaim{Lemma 4.3} Let $G$ be an abelian topological group, $S$
dense in $G$, and
      assume that either $f  \in  C_b(G,X)$ and $c_0 \not \st X$, or that $f
      \in
      C_{b,wscp}(G,X) $.
      Then to $\{t\}  \st  G$ and $\e > 0$ there exists $\delta = \delta_{t,\e}> 0$,
      $m = m_{t,\e}\in  \N$ and   $r_{1,t,\e},...,r_{m,t,\e}  \in  S$  such that

    (4.3) \qquad \qquad $\cap _{k=1}^m  E(\Delta_{r_k, t,\e} f, \delta,\{t\})  \st  E(f,\e,\{t\}) $.
\endproclaim
\demo{Proof}  We prove (4.3) case $t=0$,
 by contradiction.
Assume there is $\e_0 >0$ for which for no $\delta >0$  and
$r_1,\cdots, r_m \in S$ the (4.3) is true; take $\delta=1/2$ and
$r_1 = 0$, so $E(0,1/2,{0})$ is not contained in $E(f,\e,\{0\})$;
this means there is a $t_1 \in E(0,1/2,\{0\})$ which is not in
$E(f,\e_0,\{0\})$.  We  construct by induction a sequence
 $(\tau_n) \st S$ such that for $n\in \N$

(4.4) \qquad $||f(\tau_n) -f(0)|| > \e_0$ but

\qquad\qquad $||f(\tau_{n} + \tau_{i_k}+\cdots +\tau_{i_1}) -
f(\tau_{n+1})- f( \tau_{i_k}+ \cdots +\tau_{i_1})+ f(0)|| \le
\e_0/2^{n}$,

(4.5)\qquad  $1 \le  i_1 <\cdots < i_k < n$, $k=1, \cdots, n-1$.

\noindent Indeed, since $S$ is dense in $G$ and $f$ is continuous,
 $t_1$ above can be replaced by $\tau_1 \in S$, $\tau_1 \in
E(0,1,\{0\})$ but $\tau_1\not \in E(f,\e_0,\{0\})$. So $\tau_1$
satisfies (4.4). If $\tau_1, \cdots, \tau_{n}$ are already
constructed, we set $W_n (\e)= \cap E (\Delta _{\tau_{i_1+\cdots
+\tau_{i_k}}} f, \e, \{0\})$,
 where the intersection is over all sets of indices satisfying
(4.5). Then there is $t_{n+1} \in W_n (\e_0/2^{n+1})$ such that
$t_{n+1}\not \in E(f, \e_0, \{0\})$. By continuity of $f$, we can
replace $t_{n+1}$ by $\tau _{n+1}\in S \cap W_n (\e_0/2^{n})$ and
$\tau_{n+1} \not\in  E(f, \e_0, \{0\})$. This completes the proof
of the induction  to construct  $(\tau_n)$ satisfying (4.4),
(4.5).

\noindent Now we can proceed as in the proof of Lemma 1
(respectively Lemma 2  using [20, p. 79, IV, \S 1, (4b)] and the
Pettis-Orlicz
                      theorem [20, p.80, IV, \S 1, theorem 1]) of [4] to conclude (4.3) for the case
$t=0$ , the recurrent in [4] is ours by the remark after (2.3).
The  cases $ t\not =0$ can be shown similarly. $\square$
\enddemo
\proclaim{Lemma 4.4} Let $G$ be an abelian locally compact group,
$S$ dense in $G$, and
      assume that either $f  \in  C_b(G,X)$ and $c_0 \not \st X$, or that $f
      \in
      C_{b,wscp}(G,X) $.
      Then to compact $ K \st  G$ and $\e > 0$ there exist a compact $K_* \supset K$, $\delta = \delta_{K,\e} > 0$,
      $l = l(K,\e)\in  \N$ and   $s_1,...,s_l  \in  S$  such that

    (4.6) \qquad \qquad $\cap _{j=1}^l  E(\Delta_{s_j} f, \delta, K_*)  \st  E(f,\e,K) $.
\endproclaim
\demo{Proof} By (4.3),  to   $t\in G$, $\e >0$ there exist
$\delta_t
>0$ and $r_{1,t,\e}, \cdots, r_{m,t,\e}\in S $ such that
 $\cap _{j=1}^{m_t} E(\Delta_{r_{j, t,\e}} f, \delta_t/2,\{t\})
\st E(f,\e/2,\{t\}) $.
 Using the continuity
of $f$ and (2.3), one can find a compact neighbourhood $ W^t \st
G$ of $0$  such that  $||f(t+h) - f(t)||\le \e/2$,\, \,\,
$||\Delta_{r_{k,t,\e}} f (t+h) - \Delta_{r_{k,t,\e}} f(t)||
                        \le  \delta_t / 2$ and
 $\cap _{j=1} ^{m_t} E (\Delta _{r_{j,t,\e}} f, \delta_t/2, t+ W^t)
+W^t \st \cap _{j=1} ^{m_t} E(\Delta _{r_{j,t,\e}} f, \delta_t,
\{t\})\st E(f,\e/2,\{t\}) $ for all $h \in
                        W^t$, $1\le k\le m_t= m_{t,\e}$.
 It follows
$||f(t+\tau +h)-f(t+h)||\le ||f(t+\tau +h)-f(t)||+$
$||f(t+h)-f(t)||\le \e$,\,\,
   $h \in   W^t$,\, $\tau \in \cap _{k=1} ^{m_t} E (\Delta _{r_{k,t,\e}} f,
\delta_t/2, $ $t+ W^t)$,
 implying

(4.7) \qquad $\cap _{k=1} ^{m_t} E (\Delta _{r_{k,t,\e}} f,
\delta_t/2, t+ W^t)\st E (f, \e, t+ W^t)$,\,\,\, $t\in G$.

Let $\emptyset \not =K\st G $ be any
 compact set and $\e >0$.
 The set $\{t+ (W^t)^{\circ}: t\in K\}$ with $W^t$ satisfying (4.7) is an
 open covering for $K$. Let $ \{t_1+ (W^{t_1})^{\circ}, \cdots,  t_n+ (W^{t_n})^{\circ})\}$
be a subcovering of $K$. Set $K_*= \cup _{k=1} ^n (t_k+ W^{t_k})$
and $\delta =$ min $\{ \delta_{t_1}/2, \cdots \delta_{t_n}/2\}$.
By (4.7), one gets
 $ \cap _{j=1}^n \cap _{k=1} ^{m_{t_j}} E
(\Delta _{r_{k,t_j,\e}} f, \delta, K_*)\st \cap _{j=1}^n \cap
_{k=1} ^{m_{t_j}} E (\Delta _{r_{k,t_j,\e}} f, \delta, t_j+
W^{t_j})\st \cap _{j=1}^n E (f, \e, t_j+ W^{t_j}) \st E(f,\e,K)$.
With $l= m_{t_1}+\cdot + m_{t_n}$, one gets
    (4.6). $\square$ \enddemo

We can now generalize results of [4], from $(P_U)$ to $(L_U)$ and
from $C_{rwc}$ to $C_{b,wscp}$ in (ii) (Example 4.10):

\proclaim{Theorem 4.5} $REC(\r,X)$ satisfies $(L_U)$ and $(P_U)$
in the following
           cases :

          (i)  $U = C_b(\r,X)$  and  $c_0  \not \st X$,

          (ii)  $U = C_{b,wscp}(\r,X)$ of (4.2).
 \endproclaim
     \demo{Proof} Case  $(L_U)$: We apply Lemma 4.4 with $G = \r$ and $S =$ rationals $\Q$, getting
     (4.6). With  $\Delta_{-h} f  =  - (\Delta_h f)_{-h}$, the invariance of
     $rec \{\Delta_{1/n} f\}$ and (4.1) one gets therefore  $(\Delta_{s_1} f,...,
     \Delta_{s_l} f)  \in  REC(\r,X^l) $. (4.6) implies therefore  $f  \in  REC(\r,X)$.

 $(P_U)$ follows from $(L_U)$ with Proposition 2.12. $\square$
      \enddemo
In  the following we need:

\proclaim{Lemma 4.6}  If $K\st X$  and for each $\e>0$, there is a
relatively weakly  compact subset $K_{\e} \st X$ which is an
$\e$-net for  $K$, then $K $ is relatively weakly compact.
\endproclaim
\demo{Proof}  Let $(a_n)$ be a net,  $\st K$    and $a_n\to a^{**}
$ in the $w^*$-topology with  $a^{**}\in  X^{**}$. Let $\e
>0$ and $(b_n) \st K_{\e}$ such that

$|| a_n-b_n ||< \e$. We can assume $b_n\to b_{\e}$ weakly with
$b_{\e} \in X$. We have $|| a^{**} - b_{\e}||< {\e}$. Since $\e$
is arbitrary, $a^{**} \in X$.
\enddemo

\proclaim{Corollary 4.7}  $REC_{xy}(\r,X)$ satisfies $(L_U)$ and
$(P_U)$ for $ xy
 \in \{b, u=ub, \, rwc, rc,$
 $
  urwc, urc\}$, provided (i) or (ii) of Theorem 4.5   hold.

 These $REC_{xy}(\r,X,)$ satisfy $(\Delta)$, with
$c_0\not \st X$ if $xy=b, rwc$.
 \endproclaim
      \demo{Proof} $(L_U)$ implies $(P_U)$ by Proposition 3.9(i).

  $(L_U)$: Theorem 4.5 gives $f  \in  REC$, so  $f  \in  C_b$. The missing
$f  \in
     C_{xy}$  follows for $xy = u $ or $urc$  by the proof  of Proposition 4.2.

   $xy = rc$ or $rwc$ : By the remark after (2.3) to $\e > 0$
there exists
      a finite $F  \st  \r$ with  $E(f,\e, \{0\}) + F  =  \r$ ;  $f(\tau+s) = \Delta_s
      f (\tau) + f(\tau)$ with  $\tau  \in   E(f,\e, \{0\})$  yields therefore  $f(\r)
      \st
     \cup_{s  \in  F} (\Delta_s f)(\r)  + \{w : ||w - f(0)|| \le \e\}$.
     If $xy=rc$, this shows that  $f(\r)$ is totally bounded, $=$ relatively compact.
     If  $xy = rwc$, $f(\r)$   relatively weakly compact  follows
     from Lemma 4.6.

$(\Delta)$: $xy\in \{u, urc, urwc\}$: Proposition 3.9 (ii).

 Let $f\in L^1_{loc} (\r,X)$ and $\Delta_h f \in
      REC_{xy} (\r,X)$, $h>0$. Then with  $F= f-M_k f$,
      one has $\Delta_s F = \Delta_s f + M_k (\Delta_s f) \in  REC_{xy}
      (\r,X)$, $s >0$ by  Proposition 3.9 (i). $(L_U)$ implies
      then
       $ F \in REC_{xy}
      (\r,X)$, provided  (i) or (ii) of Theorem 4.5 hold.

(i) applies if $F\in C_b$ which is the case for $xy=b$ or $rwc$ by
[9, Example 3.5].

(ii) applies for $xy=rc$ with [10, Proposition 2.2(ii)].
      $\square$
\enddemo
    \proclaim{ Corollary 4.8} The $REC_{xy}$ of Corollary 4.7 satisfy $(L_U)$ and
    $(P_U)$
               with  $U = C_{rwc}$  or  $C_{rc} $.
\endproclaim

 \proclaim{ Corollary 4.9} If $ \emptyset\not = V  \st  REC_{xy}(\r,X)$ with (3.2) and $xy$ is as in Definition 3.2,
     then  $rec_{xy} V$ satisfies $(L_U)$ and  $(P_U)$  if (i) or (ii) of
     Theorem 4.5 hold. Except $rec V$, these $rec_{xy} V$ satisfy $(\Delta)$, with  $c_0 \not \st X$ if
       $xy = b$ or $rwc $.
\endproclaim
\demo{Proof} If  $f  \in  U$  with  all $\Delta_h f  \in  \A : =
rec_{xy} V$, then $f  \in  REC_{xy}$  by Theorem 4.5 resp.
Corollary 4.7, one has (4.6) with $\Delta_{s_j} f \in  \A $. So
$f_1,...,f_{n-1} \in \A$ implies $(f_1,...,f_{n-1},\Delta_{s_1} f,
\cdots, \Delta_{s_l} f)  \in REC(\r,X^{n-1+l})$; (4.6) implies
$(f_1,...,f_{n-1},f) \in REC(\r,X^n)$ and so  $f  \in  \A$ by
maximality.

$(\Delta)$: Cases $u,\, urwc,\, urc$ follow with Proposition 3.3
and [9, Proposition 3.1].

\noindent   The other cases follow as in the proof of Corollary
4.7.
 $\square$
\enddemo

\proclaim{Example 4.10} Even $(P_U)$ becomes false for $REC(R,X)$
and $U = C_b$
    if  $c_0   \st  X$.
\endproclaim
   \demo{Proof}    By  [1,  p. 53, (1.2)] there is an $f  \in  AP (\r,
      c_0(\N,\r))$ with $P f$ bounded but not ap; this $P f$ is  not even
      recurrent. Indeed,  $P f$ is not recurrent follows since $E(Pf, \e,
      \{0\})$
is not  relatively dense;  for otherwise $P f$ is relatively
compact by [4, Lemma 3], which would imply that $P f$ is ap [4,
Corollary 1]. $\square$
\enddemo

\proclaim{Example 4.11} There are $X$ and $M\st X$ such that $M$
is countable, bounded, relatively weakly sequentially complete and
  the closed linear hull of $M$  contains
   $c_0(\N,\r)$, but $M$ is  not  relatively weakly compact in $X$.
\endproclaim
 \demo{Proof} Take $X = \{f  \in  C(I,\cc) :
f(-\pi)=f(\pi)\}$ with $||\cdot ||_{\infty}$, $I=[-\pi,\pi]$,
   $M = \{e^{i\,n\,t}|\,I  : \, n  \in  \N \}$. Then $M$
   satisfies all the above noting that
   with the Riemann Lemma  one can show  that $M$ contains no infinite weak Cauchy
   sequence. $\square$
\enddemo

\proclaim{Example 4.12} For $ \A : = REC(\r_d,X)  \cap  C(\r,X)$,
one has :

         (a)  $\A$  satisfies $(L_U)$  with $U$, $X$ as in (i) or (ii) of Theorem 4.5.

         (b)  Not even the classical Bohl-Bohr analogue holds, i.e. already
                for $X = \cc$  the $\A$ does not satisfy $(P_b)$, then $REC(\r_d,\cc)  \cap  C(\r,\cc)
                =$
$MIN(\r_d,\cc)
         \cap  C(\r,\cc)$ by Theorem 2.3.

         (c)   $\A  \st  \m\A$ is false, also for $A_b$, already for $X = \cc$.

      \noindent This A is a first example with $(L_U)$ but not $(P_U)$; (b) gives a first
       example where the BB(AK) theorem is false for a space of generalized almost periodic functions.
\endproclaim
       \demo{Proof} (a): (4.3) of Lemma 4.3
and (4.1), with $\r_d$ instead of $\r$.

 (b) : For $g = e^{i t^2}$, one has $g  \in \A$ by Example 2.9 ; since $(Pg)(\pm \infty)$ exist,
    $O_d(Pg)$ is not minimal though Pg is bounded.

  (c) $(L_U)$ implies $(P_U)$ if  $f  \in  \A$ , $h > 0$ imply  $h M_h f
  \in
\A $. $\square$
\enddemo

  \head {\S 5  Recurrent solutions of difference-differential
  systems}\endhead

We consider systems of the form

(5.1) \qquad  $\sum^n_{k=0} \sum^m_{j=1} a_{jk}\, y^{(k)}(t+t_j) =
f(t)$,

 \noindent where  $m, r  \in  \N$, $n \in \N_0$, $ t_1<t_2<\cdots<t_m$ are real, $a_{jk}$ are complex
$r\times r$-matrices,  $f : \r\to X^r$.

\smallskip

\noindent $ y : \r\to  X^r$ is called a $solution\,\, of$ (5.1) if
$y \in$ Sobolev space $W^{1,n}_{loc} (\r,X^r)$ ([11, (2.5)]) and
satisfies (5.1) a.e  on $\r$ ($y, f$ column vectors).

 About (5.1) we have to assume ($\gamma_{\omega}(t) : = e^{i \omega t},
\omega, t \in \r$)

(5.2) \qquad det\,\,$ \sum^m_{j=1} a_{jn} \gamma_{t_j} \not\equiv
0$ on $\r$.

\smallskip

$\f = \f(\r,\cc)$ is Schwartz's space of rapidly decreasing test
functions [38, p. 146].

\proclaim{Theorem 5.1}  Assume $y$ is a solution of (5.1) on $\r$
with (5.2) as above,
      $k \in \N_0$,  $\gamma_{\omega} U_{ub}  \st  U  \st  L^ 1_{loc}(\r,X)$ for all $\omega
     \in \r$,  $U_{ub} * \f  \st  U$, $U$ such that for $g \in REC_{urc}(\r,X^r)$  the
       $rec_{urc}\, \{g_1,...,g_r\}$ satisfies $(L_U)$ (Definitions 3.2/4.1).

    \noindent If then  $ f \in \m^{k+1} REC_{urc}(\r,X^r)$  and  $y \in \m^k(
(C_{ub}
       (\r,X^r))  \cap  U^r)$,  then  $y \in \m^k (REC_{urc}(\r,X^r))$.
\endproclaim
\demo {Proof} Case k=0 : By assumption, $g = (g_1,....,g_r) : =
M_1(f) \in REC_{urc}(\r,Y)$, $Y : = X^r $. Then $V : =
\{g_1,...,g_r\}$ satisfies (3.2), so  $\A : =  rec_{urc} V $ is
well defined and a $\lambda$-class by Corollary 3.4, $V \st
REC_{urc}(\r,X)$. With suitable $\chi \in L^1(\r,\cc)$ with
compact support one has    $f_j
* \chi = M_1 f_j = g_j \in \A$   for $1 \le j\le r $.

 Now the proof for Lemma 2.4
  of [7], especially p. 126, works
also, instead of the (i),, (ii) or (iii) there, with the more
general $(L_U)$ for $\A$, with $U$ as in Theorem 5.1; this implies
that Theorem 2.5 and Corollary 2.6 of [7] hold also with this
$(L_U)$ for $\A$ instead of the (i), (ii) or (iii), and $\Omega_j
\in U$.

  Now all assumptions of this generalized Corollary 2.6 of [7] are
fulfilled, one gets  $y_j \in rec_{urc} V$ for $1 \le j \le r$,
with (3.2) for $rec_{urc} V$  this gives  $y \in REC_{urc}(\r,Y)
$.

 The case $k
> 0$ follows by applying the above to $z : = M_h y$, $h =
(h_1,...,h_k)$ with  $h_j
> 0 $.  $\square$
\enddemo

 \head {Special cases}\endhead

 \noindent (a) \,\,\,  $U =$ totally ergodic functions $\T\E$, (for $\E$ see after Proposition 4.2)

  (5.3)\qquad
$\T\E(\r,X)  : =  \{f \in \E(\r,X) : \gamma_{\omega} f \in
\E(\r,X)$ for all real $\omega \}$ :

\noindent  $\A = rec_{urc}\, \{g_1,...,g_r\}$ of Theorem 5.1
satisfies even $(L_{\E})$ by Proposition 4.2;
  $\E_{ub} * \f  \st  \E$  by [7, Lemma 2.2 (A)], $\gamma(f*\va) = (\gamma f)*
  (\gamma \va)$.

\noindent (b) \,\,\,  $U = C_b(\r,X)$ with $ c_0  \not \st  X $:
Corollary 4.9 gives $(L_b) $.

\noindent (c) \,\,\, $U = C_{rwc}(\r,X)$: By Corollary 4.9 one has
$(L_U)$; $U_{ub}$ is even  a $\lambda$-class by [7, p. 119],
$U_{ub} * \f \st  U$  follows with
   [7, Lemma 2.2 (A)].

\proclaim {Corollary 5.2}  If $X = \cc$, $y$ a solution of (5.1)
with (5.2), $y \in
               C_{ub}(\r,\cc^r)$ and  $f  \in  REC_b(\r,\cc^r)$, then  $y
               \in
                REC(\r,\cc^r) $.
\endproclaim
  \noindent  Here we use $REC_b(\r,X^r)  \st  \m REC_{ub}(\r,X^r)$ , $REC_{urc}(\r,\cc^r)
     = REC_{ub}(\r,\cc^r)$.

\proclaim {Corollary 5.3}  If  to  $F  \in  C_{ub}(\r,\cc)$  there
exists one  $h_0 > 0$ with
             difference  $F_{h_0}  -  F$  recurrent, then $F$ itself is recurrent.
\endproclaim

\proclaim {Corollary 5.4} If $y$ is a distribution-solution of
(5.1) with (5.2), $U$ and
       $(L_U)$ as in Theorem 5.1,  $f  \in  \h'_{REC_{urc}(\r,X^r)}$, with $y_j
        \in  \h'_U  \cap  \h'_{L^{\infty}(\r,X)} $ if $1 \le j \le r$, then  $y  \in  \h'_{REC_{urc}(\r,
        X^r)} $.
\endproclaim
  Here for  $\A  \st  L^1_{loc}(\r,Y)$ the  $\h'_{\A} $ contains all
distributions $T  \in  \h'(\r,Y)$  with convolution  $T*\va  \in
\A$ for $\va$ test function  $\in
 \h(\r,\cc)$ (see [11]);  one has  $\A  \st  \h'_{\A}$  for  $\A = REC(\r,X)$, $REC_b, REC_u,
      REC_{rc}$ and combinations (for $REC_{rc}$ e.g. as in [10, Proposition 2.2 (i)]).

\demo{Proof} For   $\va \i  \h(\r,\cc)$, Theorem 5.1 can be
applied to  $y*\va$ with $k = 0$, using  $\h'_{L^{\infty}} =
\h'_{C_{ub}}$ of [11, Proposition  2.9] and Proposition  3.6(i).
$\square$

\enddemo

\proclaim{Remarks 5.5}

 \noindent (i) $y'+y=0$ shows that Theorem 5.1 becomes
false for
    half-lines; $(rec_{urc} \, V)|[0,\infty)$ is no longer a $\lambda$-class [7, p.
    117, (1.III)]. But see Proposition 5.9 below.

 \noindent (ii)  If  $| det\,(\sum^m_{j=1} a_{j,n} \gamma_{t_j})| \ge
\delta_0
> 0$ on $\r$,
     then the assumption on $y$ can be weakened to  $y  \in  \m^k(C_b(\r,X^r)
      \cap  U^r)$,  i.e. the "uniform continuity" can be dropped (special case :
      m=1) : [7, p. 131 (e)] ; but see [7, p. 129, Example 2.8(ii)].

 \noindent (iii)  $f  \in  \m^{k+1} REC_{urc}$  can be weakened to: There
exists one
      $\psi_0  \in  L^1(\r,\cc)$ with compact  support such that  $f*\psi_0
      \in
      \m^k(REC_{urc}(\r,X^r)$.

 \noindent (iv) If  $r > 1$,  $f  \in  \m^{k+1}(REC_{urc}(\r,X))^r$ is not
enough in Theorem
       5.1  because of (3.2).

 \noindent
  (v)  Theorem 5.1 can be extended to
integro-difference-differential systems as in  [7, Corollary 2.6].
\endproclaim

The classical Esclangon-Landau theorem which says that bounded
solutions of linear differential equations with bounded right-hand
side have bounded derivatives, now can be extended to recurrent
solutions of (5.1) :

\proclaim{Proposition 5.6} If $y$ is a solution of (5.1) with
(5.2), $n\ge 1$, $a_{j,n} =0$ for  $j \not =j_0$, with one
$t_{j_0} = 0$, $f  \in  \m C_u (\r,X^r)$, there are $s$, $s^* \in
\N_0$
     with  $y \in (\m^s C_u
      (\r,X^r)) \cap \m^{s^*} REC(\r,X^r)$, then $y^{(k)} \in REC(\r,X^r)$ for
      $0\le k < n $, $y^{(n)}  \in  \m REC(\r,X^r)$.
\endproclaim
  Here  $g$  bounded or $S^p$-bounded  respectively ergodic implies  $g
\in \m C_u$
     respectively $\m^2 C_u $.

\demo{Proof} As for Theorem 5.9 of [9]; instead of Lemma 2.2 of
[11] used there, which cannot be applied to $\A = REC$, one can
use Lemma 5.7 below.  Proposition 2.9 of [9] can be here applied,
since all REC-
            spaces are uniformly closed.   $\square$
\enddemo

\proclaim{Lemma 5.7} If  $\A  \st  L^1_{loc}(\jj,X)$  satisfies

(5.4)\qquad  $f \in \A$ , $h > 0 $ imply  $(\Delta_h f)/h  \in
\A$,

\noindent then for $n \in \N$, $F  \in  \A \cap
W^{1,n}_{loc}(\jj,X)$ one has $F^{(n)} \in \m^n \A$ .
\endproclaim
\demo{Proof}  As for Lemma 2.2(c) of [11], with $W^{1,n}_{loc}$ of
[11, (2.5)]. $\square$
\enddemo
\proclaim{Remarks 5.8} (i) In Proposition 5.6 the REC can be
replaced everywhere by
     $ V$, with  $V \in \{REC_b,\, REC_u, \,REC_{ub},
     \,REC_{rc},\,
     REC_{urc}\}$.

(ii)   Also extensions to solutions on half lines are possible:
(see Proposition 5.9 and the proof of Corollary 5.10).

\endproclaim

  For half lines  $\jj = [\al,\infty) \not =\r $ we can only treat a
special case of (5.1) (see Remark 5.5(i)) :

(5.5) \qquad $ Ly : = y^{(n)} + \sum^{n-1}_{k=0} a_k y^{(k)}  =
f$,

\noindent with $n, r  \in \N$, $a_k$ complex constant $r \times r$
- matrices, solution on $\jj$ defined as after (5.1);

  $spectrum$\,\,\,  $\sigma(L) : = \{ \lambda  \in  \cc :$ det $(\lambda^n  +
                                                  \sum^{n-1}_{k=0} \lambda^k) = 0\}$,

(5.6)  \qquad         $REC(\jj,X)  : =  REC(\r,X)|\jj$.

  For $\A, U  \st  X^{\jj}$, $\A$ $\,satisfies\, $ $(P'_U)$  means (see Definition 4.1, [8,
Proposition 3.12])

 $(P'_U)$ $\, f  \in  \m\A$ and $Pf  \in U$  imply  $Pf \in
A $.

(5.7) \qquad  $O_q(\jj,X) : = \{g  \in L^1_{loc}(\jj,X) :$ sup
$_{t \in \jj} ||g(t)/(1+|t|)^q|| < \infty \} $.

\proclaim{Proposition 5.9} Assume $\emptyset \not  = V  \st
REC_{urc}(\r,X)$, $V$ satisfies (3.2),
     and define $\A = (rec_{urc} V)|\jj$; assume that $\A$ satisfies $(P'_U)$ with
     constants $X  \st  U  \st  L^1_{loc}(\jj,X)$, $U$ linear with  $(\Gamma)$ of
     Definition
     3.1.

   If then  $\sigma (L)  \st  \{\lambda  \in\cc :$  Re $\lambda \ge 0 \}$,  $y$ is
a solution
     of (5.5)  on $\jj$ with  $f  \in  \m(\A^r)$,  $y  \in  \m^s O_q(\jj,X^r)$   for some
      $s, q  \in  N_0$  and  $y^{(k)}  \in  U^r $ for $ 0 \le k \le n-1$,  then  $y  \in  \A^r$,
      $\st  REC_{urc}(\jj,X^r) $.
\endproclaim
\demo {Proof} This follows from Theorem 5.10 of [8], all the
assumptions there are here fulfilled :

  $\A$ is linear, positive-invariant, with $(\Gamma)$, $(\Delta)$, $X  \st
  \A \st \m\A$
 by Corollary 3.4, Corollary 3.6 and [10, Proposition 2.3], $\A$ is  uniformly closed   by Proposition 3.3 and
[6, Proposition 2.1.7].

  So with $\A  \st   C_b$  and [8, Proposition 4.15] all assumptions of
(a) (of Theorem 5.1 of [8]) are fulfilled.  Since A has $(\Gamma)$
and $(\Delta)$, by [8, p. 678 [Proposition 4.1(ii)]] also $\m\A$
satisfies $(\Gamma)$, therefore all assumptions of (b) of [8,
Theorem 5.10] hold, one gets indeed  $y \in  A^r$. $\square$
\enddemo

\proclaim{Corollary 5.10}  If $y$ is a solution of (5.5) on $\jj$
with $\sigma (L) \st \{$ Re
      $\lambda \ge 0\}$,  $f = F|\, \jj \in  REC_{urc}(\jj,X^r)$, and  $y  \in  U^r$, where  $U =
      \T\E(\jj,X)$ (see (5.3)), $C_{rwc}(\jj,X)$, or  $C_b(\jj,X)$ with  $c_0  \not \st X$,
      then  $y  \in  ((rec_{urc} \,\{F_1,...,F_r\})|\jj)^r$,   $\st  REC_{urc}(\jj,X^r) $.
\endproclaim
\demo{Proof} If $f = F|\jj$, one can use  $V : = \{F_1,...,F_r\}$
for the $\A$ of Proposition 5.9. $y  \in  \m^s O_q $ follows with
$\T\E \st  \E  \st \m C_b \st \m^2 C_u$ of  [9, (2.4)]. $y^{(k)}
\in U^r$  for $0 <k < n$ can be obtained with [9, Theorem 5.9,
Remark 5.8(b)], that $C_{rwc}$ is uniformly closed follows as in
[7, p. 119].

 $ \A$ satisfies $(P'_U)$:
With the properties of $\A$ (proof of Proposition 5.9), for  $U =
\T\E$ this follows from [8, Corollary 4.2].
 If $U = C_{rwc}(\jj,X)$, with [8, Proposition 3.12]  and $U$ linear it
is enough to show $(P_U) $; if $f \in \A$   with  $Pf  \in  U$,
then $f = F|\jj$ with  $F \in rec_{urc} V $. Now

(5.8) \qquad $(P_{\beta} F)(\r)  \st$ norm closure of $((P_{\beta}
F)([\gamma,\infty ))
                            -   (P_{\beta} F)([\gamma,\infty)))$

 \noindent  if  $\beta$, $\gamma  \in \r$   and  $F  \in  REC(\r,X)$,  $(P_{\beta} F)(t) : =
\int^t_{\beta} F(s)\, ds $ (proof as in [8, Lemma 4.4]).

\noindent With this one gets here   $PF  \in C_{rwc}(\r,X)$ with
[8, Proof of Proposition 4.7] ;  Corollary 4.9 gives $PF \in
rec_{urc} V $ and then $Pf \in \A$.

 $ U = C_b$ follows with Corollary 4.9 in the same way or with
[8, Lemma 4.4].   $\square$
\enddemo

\proclaim{Remarks 5.11} (i) Proposition\,5.9 and Corollary 5.10
become false without

    $\sigma(L)  \st   \{$ Re $\lambda  \ge  0 \}$.

(ii) In Proposition\, 5.9/Corollary 5.10 no explicit assumption "
$y$ uniformly
      continuous" is needed, in accordance with Remark 5.6(ii).

(iii)  $y  \in  (rec \{f_1,...,f_r\} | \jj )^r $  an be
interpreted as " $y$ is recurrent
      as $f$" .

(iv)  Also the other results of [8, \S 5] can be extended to
recurrent
     solutions; for example Theorem 5.13/Corollary 5.15 of [8] would
     strengthen Proposition 5.9/Corollary 5.10 above for  $r = 1$.
\endproclaim

\head {\S 6     Open Questions}\endhead

1. Is $(\Gamma)$ true for $REC_{rc}$, $REC_{urwc}$, $REC$ or
$rec_{urc}$, $rec_{urwc}$, $rec$?

2. Is Theorem 4.5 true for $G = \r^2 $?

3. Is Theorem 5.1 true for $REC_{xy}$, $xy = urwc$ or $ub$,
instead of
                        $REC_{urc} $?

4. Is Corollary 5.3 true if only $f  \in  C_u$ or $C_b $?

5. Is  $n = 2$  enough in (3.2)?

6. Is  $REC(\r,X)|[0,\infty)  =  REC([0,\infty),X)$ where for the
latter in
      Definition 2.1

      \qquad the G is replaced by $[0,\infty) $?

7. Do there exist $f  \in  REC$ with  $f(\cdot) + f(-\cdot)  \not
\in REC $?

 8. Is  $f+g  \in  REC_{xy}$  equivalent with  $(f,g)  \in   REC_{xy}$   for  $f
 \in
          REC_{xy}$, $g =  \gamma_{\omega}$

          \qquad (or  $\in  AP$,  $\in  REC_{xy}$)?

9. Do $REC$ or $rec $ satisfy $(\Delta)$, is $c_0 \not \st X$
superfluous in Corollary 4.7 or 4.9?

10. Does there exist   $f  \in  REC$ with $Pf \not \in C_{rwc}$,
but $ Pf\in
     C_{b,wscp}$ of (4.2) and

     \qquad $c_0  \st$  smallest $X$ containing $Pf(\r) $?

11. For what $xy$  exist  $f  \in  REC_{xy}$ with $O(f)  \not \st
REC $?

12. Is $C_{b,wscp}$ of (4.2) linear, a $\lambda$-class (Definition
3.1)?

13. For what $xy$ has one   $\m^n \A  \st  \h'_{\A}$ , $\A =
REC_{xy} $?
    ($\h'_{\A}$: after Corollary 5.4.)

14 What is $rec_b \{ e^{i\,t^2}\}$, for  $\r_d$? (See Example
2.9.)

\indent School of Math. Sci., P.O. Box No. 28M, Monash University,
 Vic. 3800.

\indent E-mail "bolis.basit\@sci.monash.edu.au".

\indent Math. Seminar der  Univ. Kiel, Ludewig-Meyn-Str., 24098
Kiel, Deutschland.

\indent E-mail "guenzler\@math.uni-kiel.de".

\enddocument